%% file: chenrubin.tex
\documentclass[11pt]{article}

 \usepackage{a4}                          
 \usepackage{amsmath}
 \usepackage{amssymb}
 \usepackage{graphicx}
 \usepackage{verbatim}
  
 \newtheorem{thm}{Theorem}[section]

 \newtheorem{prop}[thm]{Proposition}

 \newtheorem{lemma}[thm]{Lemma}
 \newtheorem{rem}[thm]{Remark}

 \title{The Chen-Rubin conjecture in a continuous setting}
 \author{Christian Berg and Henrik L. Pedersen\footnote{Research supported by the
     Carlsberg Foundation}}

 \date{\today}

\begin{document}

 \maketitle

 \begin{abstract} We study the median $m(x)$ in the gamma distribution
   with parameter $x$ and show that $0<m'(x)<1$ for all
   $x>0$. This proves a generalization of a conjecture of Chen and
   Rubin from 1986: The sequence $m(n)-n$ decreases for $n\geq 1$. We also describe the
   asymptotic behaviour of $m$ and $m'$ at zero and at infinity.
 \end{abstract}



\section{Introduction}
The gamma distribution with parameter $x$ has density with respect to
Lebesgue measure on $(0, \infty )$ given by $e^{-t}t^{x-1}/\Gamma (x)$.
We consider the median $m(x)$ of this distribution which is defined implicitly as 
$$
\int_0^{m(x)}\frac{e^{-t}t^{x-1}}{\Gamma (x)}\, dt=\frac{1}{2},
$$
or
\begin{equation}
\label{eq:def_m_0}
\int_0^{m(x)}e^{-t}t^{x-1}\, dt=\frac{1}{2} \int_0^{\infty}e^{-t}t^{x-1}\, dt.
\end{equation}
We of course also have 
\begin{equation}
\label{eq:def_m_infinity}
\int_{m(x)}^{\infty}e^{-t}t^{x-1}\, dt=\frac{1}{2} \int_0^{\infty}e^{-t}t^{x-1}\, dt.
\end{equation}
We show that $m$ is continuous and increasing. This is a consequence
of a result about general convolution semigroups of probabilities on
the positive half-line, that we give in Section \ref{sec:semigroup}. There we
also show that $m$ is real analytic and that $m$ satisfies a
certain differential equation.

We shall mainly study $m$ through the function 
\begin{equation}
\varphi(x)\equiv \log\frac{x}{m(x)}, \quad{x>0}.
\end{equation}
This function appears if we make the substitution $u=\log (x/t)$ in
the relation (\ref{eq:def_m_infinity}). We get:
\begin{equation}
\label{eq:def_phi_1}
\int_{-\infty}^{\varphi (x)}e^{-x(e^{-u}+u)}\, du =\frac{1}{2} \int_{-\infty}^{\infty}e^{-x(e^{-u}+u)}\, du.
\end{equation}
Chen and Rubin (see \cite{CR}) studied the median of the gamma
distribution and proved that $x-1/3<m(x)<x$ for $x>0$.
The relation (\ref{eq:def_phi_1}) was also used in \cite{CR} to establish that
$\varphi(x)>0$, or equivalently $m(x)<x$. (It follows by observing
that $\int_{-\infty}^{0}e^{-x(e^{-u}+u)}\,
du<\int_{0}^{\infty}e^{-x(e^{-u}+u)}\, du$, which is true because
$\sinh u>u$ for $u>0$.)

Chen and Rubin furthermore conjectured that the sequence $m(n)-n$
decreases. This conjecture has recently been verified by Alzer (see
\cite{A}). He proved that $m(n+1)-\alpha n$ decreases for all $n\geq 0$ exactly when
$\alpha \geq 1$ and increases exactly when $\alpha \leq m(2)-\log 2$. Note that $m(1)=\log 2$.

In this paper we investigate the properties of $m$ as a function on
$(0, \infty)$. Chen and Rubin's conjecture follows from the relation
$$
m'(x)<1 \quad \text{for all}\quad x>0
$$
and it is our main goal to verify this
relation. (See Theorem \ref{thm:cr_phi}.) 

\begin{figure}[htb]
\begin{center}
\includegraphics[scale=0.6]{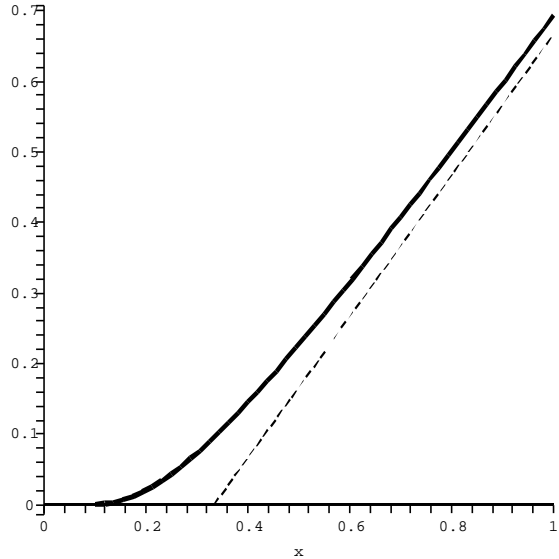}
\end{center}
\caption{The graph of $m$}
\label{figure:m}
\end{figure}

We have
drawn the graph of
$m$ in Figure \ref{figure:m}. (The dotted line is given by $x-1/3$.)

In terms of $\varphi$ Chen and Rubins conjecture takes the form
$1-x\varphi'(x)<e^{\varphi(x)}$. As a crucial step towards this result
we show that $x\varphi(x)$ decreases and that
$$
\frac{1}{3}<x\varphi(x)<\log 2
$$
for all $x>0$ (see Proposition \ref{prop:xphi}). In terms of $m$ this relation can be rewritten as 
$$
xe^{-\log 2/x}<m(x) < xe^{-1/3x}.
$$
If we use that $e^{-a}<1-a+a^2/2$ for $a>0$ then we see that in fact
$$
m(x)<x\left( 1-\frac{1}{3x}+\frac{1}{18x^2}\right) =
x-\frac{1}{3}+\frac{1}{18x}
$$
for all $x>0$. This result improves a result of Choi about $m(n+1)$, see
\cite[Theorem 1]{choi}, even though it is there claimed to be best possible.

The asymptotic behaviour of $m(x)$ for $x\to 0$ and for $x\to \infty$ is given in Section
\ref{sec:asymptotic}. From this we can deduce that $m^{(k)}(x)>0$ for
$x$ close to 0 for each $k$ while $(-1)^km^{(k)}(x)>0$ for
$x$ sufficiently large and $k\geq 2$. It is reasonable to believe that
$m''(x)>0$ for all $x>0$, i.e.\ $m$ is convex, but the higher derivatives of odd order change
sign.

We also relate our work to papers of Ramanujan, Watson and others on Ramanujan's
rational approximation to $e^x$, see Section \ref{sec:ramanujan}.

\section{Medians of convolution semigroups on the half-line}
\label{sec:semigroup}
A family $\{ \mu_x\}_{x>0}$ of probabilities concentrated on
$[0,\infty)$ is called a convolution semigroup if it has the properties
\begin{enumerate}
\item[(i)] 
$\mu_{x}([0,\infty))= 1$ for all $x>0$;
\item[(ii)] 
$\mu_x\ast \mu_y=\mu_{x+y}$ for all $x, y>0$; 
\item[(iii)] 
$\mu_x \to \delta_0$ for $x\to 0$ in the vague topology (here $\delta_0$ denotes the Dirac
mass at zero).
\end{enumerate}
The gamma distributions $\{ e^{-t}t^{x-1}/\Gamma (x)\, dt\} _{x>0}$ is an
example of such a convolution semigroup.

A probability measure $\mu$ on $[0,\infty)$ has median $m$ if 
$$
\mu ([0,m])=\frac{1}{2}.
$$
Of course a probability measure may not have a median and if it exists
it may not be unique. However, if the
measure has density w.r.t.\ Lebesgue measure on $[0, \infty)$ then the
median exists: this follows from the fact that $M\to \mu ([0,M])$ is
continuous and increases from 0 to 1. If the density is strictly
positive almost everywhere on $[0,\infty)$ then the median is unique.
\begin{prop}
\label{prop:semigroup}
Let $\{\mu_x\}_{x>0}$ be a convolution semigroup of probabilities on
$[0,\infty)$ having a.e.\ strictly positive 
densities w.r.t.\ Lebesgue measure on $[0, \infty)$. Then the median
$m(x)$ of $\mu_x$ is a continuous and strictly increasing function on $(0,
\infty )$.
\end{prop}
{\it Proof.} We write $d\mu_x(t)=g_x(t)\, dt$. As we have seen above, $m$ exists as a function on the
positive half-line. We now show that $m$ is increasing. Let $x>0$ and
$h>0$. Then 
\begin{eqnarray*}
\int_0^{m(x)} g_{x+h}(t)\, dt & = & \int_0^{m(x)} g_{x}\ast g_{h}(t)\,
dt\\
& = & \int_0^{m(x)} \int_0^{t}g_{h}(t-s)g_{x}(s)\, ds \, dt \\
& = & \int_0^{m(x)} \int_s^{m(x)}g_{h}(t-s)\, dt \, g_{x}(s)\, ds \\
& = & \int_0^{m(x)} \int_0^{m(x)-s}g_{h}(t)\, dt \, g_{x}(s)\, ds \\
& < & \int_0^{m(x)} \int_0^{\infty}g_{h}(t)\, dt \, g_{x}(s)\, ds \\
& = & \int_0^{m(x)} g_{x}(s)\, ds \quad = \quad \frac{1}{2}.
\end{eqnarray*}
Since $m(x+h)$ by definition satisfies 
$$
\int_0^{m(x+h)} g_{x+h}(t)\, dt =\frac{1}{2},
$$
we must therefore have $m(x)<m(x+h)$. This shows that $m$ is
strictly increasing.

Concerning the continuity we first notice that for any $A>0$, 
$$
\int_0^{A} g_{x_0+h}(t)\, dt \to \int_0^{A} g_{x_0}(t)\, dt 
$$
as $h\to 0$. In fact, using the same computation as above we get
\begin{eqnarray*}
\int_0^{A} g_{x_0+h}(t)\, dt & = & \int_0^{A} g_{x_0}\ast g_{h}(t)\,
dt\\
& = & \int_0^{A} \int_0^{A-s}g_{h}(t)\, dt \, g_{x_0}(s)\, ds \\
& \to & \int_0^{A} \int_0^{A-s}\, d\delta_0(t) \, g_{x_0}(s)\, ds \\
& = & \int_0^{A} g_{x_0}(s)\, ds,
\end{eqnarray*}
as $h\to 0$, because of the vague convergence. (All measures have the
same total mass, so vague convergence implies weak convergence.) 
Now let $x_0>0$ and let $\epsilon>0$ be given. For
$A=m(x_0)+\epsilon$ we have, as $h\to 0_+$,
$$
\int_0^{m(x_0)+\epsilon} g_{x_0+h}(t)\, dt \to
\int_0^{m(x_0)+\epsilon} g_{x_0}(t)\, dt > \frac{1}{2},
$$
so that $m(x_0+h)\leq m(x_0)+\epsilon$ for all sufficiently small
and positive $h$. Similarly we find $m(x_0-h)\geq m(x_0)-\epsilon$, and
this shows that $m$ is continuous. \hfill $\square$

We now specialize to consider the medians $m(x)$ of the gamma
distributions. We notice
\begin{prop}
\label{prop:m_differentiable}
The median $m(x)$ of the gamma distributions is real analytic.
\end{prop}
{\it Proof.} We consider the $C^1$-function
$$
F(x,y)=\int_0^{y}e^{-t}t^{x-1}\frac{1}{\Gamma(x)}\, dt, \quad x,y>0.
$$
The median is implicitly defined as $F(x,m(x))=1/2$. The fact that the
continuous function $m$
is $C^1$ follows from the implicit
function theorem, which yields the following differential equation for
$m$:
\begin{equation}
\label{eq:m_differential_eq}
e^{-m(x)}m(x)^{x-1}m'(x) = \frac{1}{2}
  \Gamma'(x) -\int_0^{m(x)}(\log t)\, e^{-t}t^{x-1}\, dt,
\end{equation}
which shows that $m$ satisfies a differential equation of the form 
$$
m'(x)=G(x,m(x)),
$$
with $G(x,y)$ being real analytic for $x,y>0$. Therefore $m$ is real
analytic. \hfill $\square$

\begin{rem}
From (\ref{eq:m_differential_eq}) it seems difficult to deduce
monotonicity properties of $m$; e.g.\ it is not at all clear that
$m'(x)>0$, which we know from Proposition \ref{prop:semigroup}.
For another derivation of $m'(x)>0$ see Proposition \ref{prop:xphi}.
\end{rem}

\section{Uniform results}
\label{sec:uniform}
In this section we prove the generalized version of Chen and Rubin's
conjecture (Theorem \ref{thm:cr_phi}). 

Proposition \ref{prop:key} below is the key to our results. Before
stating it we need some notation. 

\begin{figure}[htb]
\begin{center}
\includegraphics[scale=0.6]{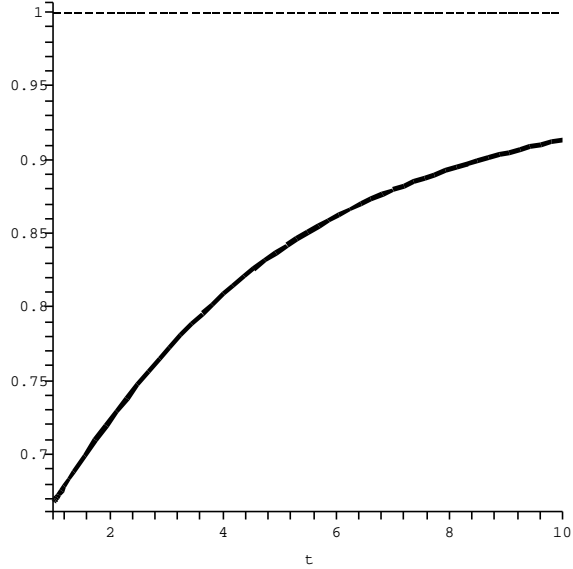}
\end{center}
\caption{The graph of $\xi$}
\label{figure:xi}
\end{figure}

We consider the function $f(x)=e^{-x}+x$. It is easily seen
that $f(x)$ decreases for $x<0$ and then increases for $x>0$. Hence
$f$ has an inverse on $(-\infty, 0]$, which we call $u$, and an inverse on
$[0,\infty)$, which we call $v$.
The function $u$ is defined for $t\geq 1$ as $u(t)\leq 0$ and 
$$
e^{-u(t)}+u(t)=t, \quad t\geq 1;
$$
the function $v$ as $v(t)\geq 0$ and 
$$
e^{-v(t)}+v(t)=t, \quad t\geq 1.
$$
The following function $\xi$ plays an important role:
\begin{equation}
\label{eq:xi}
\xi(t)=\frac{1}{1-e^{-v(t)}}+\frac{1}{1-e^{-u(t)}}= u'(t)+v'(t).
\end{equation}
The function $\xi(t)$ is defined for $t>1$ by (\ref{eq:xi}), but for
$t=1$ this expression yields $\infty -\infty$ and a closer study is
necessary. The following result holds and will be proved in Section
\ref{sec:xi}.
\begin{prop}
\label{prop:xi_properties}
The function $\xi(t)$ defined for $t>1$ by (\ref{eq:xi}) has a
holomorphic extension to the cut plane $\mathbb C \setminus \{ x\pm
2\pi i\, | \, x\geq 1\}$. The following holds:
\begin{enumerate}
\item[(i)]
$\lim_{t\to -\infty}\xi (t)=0$, $\lim_{t\to \infty}\xi (t)=1$, $\xi(t)$
is increasing;
\item[(ii)]
$\lim_{t\to \infty}\xi^{(n)} (t)=0$ for $n\geq 1$;
\item[(iii)]
$\xi(t)$ is concave for $t\in [1,\infty)$;
\item[(iv)]
$\xi(1)=2/3$, $\xi'(1)=8/135$ and $\xi''(1)=-16/2835$.
\item[(v)] All points on  $\{ x\pm
2\pi i\, | \, x\geq 1\}$ are singular points for $\xi$ and the series
$$
\sum_{k=0}^{\infty}\frac{\xi ^{(k)}(1)}{x^{k+1}}
$$
diverges for all $x$.
\end{enumerate}
\end{prop}
In Figure \ref{figure:xi} we have shown the graph of $\xi(t)$ for
$t\geq 1$.

\begin{prop}
\label{prop:key}
We have
\begin{equation}
\label{eq:key}
2\int_{0}^{\varphi (x)}e^{-x(e^{-u}+u)}\, du = \int_{1}^{\infty}\xi(t)e^{-xt}\, dt,
\end{equation}
where the function $\xi$ is defined in (\ref{eq:xi}).
\end{prop}
{\it Proof.} From (\ref{eq:def_phi_1}) we get
$$
2\int_{0}^{\varphi (x)}e^{-x(e^{-u}+u)}\, du =\int_{0}^{\infty}e^{-x(e^{-u}+u)}\, du-\int_{-\infty}^{0}e^{-x(e^{-u}+u)}\, du.
$$
In the first integral on the right hand side we make the substitution
$u=v(t)$ and in the second the substitution $u=u(t)$, where $u(t)$ and
$v(t)$ denote the functions above. In this way we get
\begin{eqnarray*}
2\int_{0}^{\varphi (x)}e^{-x(e^{-u}+u)}\, du
& = & \int_{1}^{\infty}e^{-xt}\left(
  \frac{1}{1-e^{-v(t)}}+\frac{1}{1-e^{-u(t)}} \right)\, dt\\
& = & \int_{1}^{\infty}\xi(t)e^{-xt}\, dt. 
\end{eqnarray*}
\hfill $\square$

It is easy to see that in terms of the function $\varphi$ the
relation $m'(x)<1$ takes the form
$$
1-x\varphi'(x)< e^{\varphi(x)}.
$$
\begin{thm}
\label{thm:cr_phi}
We have 
$$
1-x\varphi'(x)< e^{\varphi(x)}
$$
for all $x>0$.
\end{thm}
Before proving this result, we give the following Lemma.
\begin{lemma}
\label{lemma:key}
The relation 
$$
2\int_{0}^{x\varphi (x)}e^{-s}e^{x(1-e^{-s/x})}\, ds =\frac{2}{3}+
\int_{1}^{\infty}\xi'(t)e^{-x(t-1)}\, dt
$$
holds for all $x>0$.
\end{lemma}
{\it Proof.} In the left hand side of the relation in Proposition
\ref{prop:key} we make the substitution $s=xu$ and get in this way
$$
\frac{2}{x}\int_{0}^{x\varphi (x)}e^{-s-xe^{-s/x}}\, ds = \int_{1}^{\infty}\xi(t)e^{-xt}\, dt.
$$
On the right hand side we perform integration by parts, and this gives
\begin{eqnarray*}
2\int_{0}^{x\varphi (x)}e^{-s-xe^{-s/x}}\, ds &=&
\int_{1}^{\infty}\xi(t)xe^{-xt}\, dt\\
&=&e^{-x}\xi(1)+\int_{1}^{\infty}\xi'(t)e^{-xt}\, dt.
\end{eqnarray*}
Here $\xi(1)=2/3$ (see Proposition \ref{prop:xi_properties}) and the
assertion of the lemma now follows by multiplication by
$e^{x}$. \hfill $\square$

\begin{figure}[htb]
\begin{center}
\includegraphics[scale=0.6]{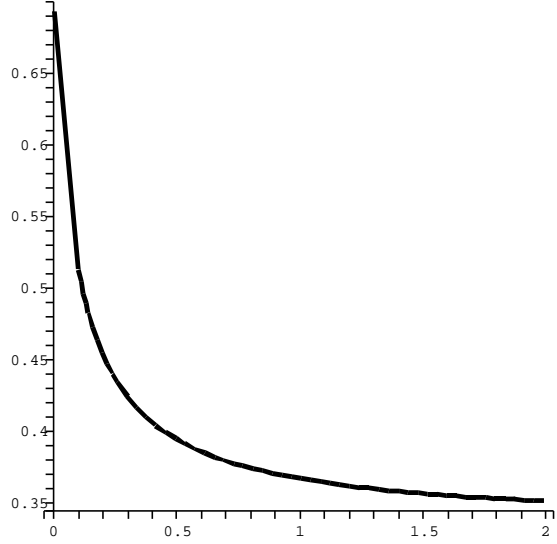}
\end{center}
\caption{The graph of $x\varphi(x)$}
\label{figure:xphi}
\end{figure}

\begin{prop}
\label{prop:xphi}
The function $x\to x\varphi(x)$ decreases for $x>0$ and 
we have
\begin{eqnarray*}
\lim_{x\to 0_+}x\varphi (x) & = & \log 2,\\
\lim_{x\to \infty}x\varphi (x) & = & \frac{1}{3}.
\end{eqnarray*}
In particular $\varphi$ is decreasing from $\infty$ to 0 (and $m$ is
increasing). (The graph of $x\varphi(x)$ is shown in Figure \ref{figure:xphi}.)
\end{prop}

{\it Proof.} From Lemma \ref{lemma:key} we find by differentiation, that
\begin{eqnarray*}
\lefteqn{e^{-x(\varphi(x)-1+e^{-\varphi(x)})}\,(x\varphi(x))' = }\\
&&-\int_{0}^{x\varphi (x)}e^{-s}e^{x(1-e^{-s/x})}\left(
  1-\left( 1+\frac{s}{x}\right) e^{-s/x}\right)\, ds \\
&&-\frac{1}{2}\int_{0}^{\infty}te^{-xt}\xi'(t+1)\, dt.
\end{eqnarray*}
Since $1-(1+a)e^{-a}>0$ for $a>0$ and $\xi'(t+1)>0$ for $t>0$ (see
Proposition \ref{prop:xi_properties}) we
have 
$$
A(x)\equiv \int_{0}^{x\varphi (x)}e^{-s}e^{x(1-e^{-s/x})}\left(
  1-\left( 1+\frac{s}{x}\right) e^{-s/x}\right)\, ds >0
$$
and
$$
B(x)\equiv \frac{1}{2}\int_{0}^{\infty}te^{-xt}\xi'(t+1)\, dt >0.
$$
Therefore we have
\begin{equation}
\label{eq:smart}
(x\varphi(x))' = -e^{x(\varphi(x)-1+e^{-\varphi(x)})}(A(x)+B(x)),
\end{equation}
which is a negative quantity, so $x\varphi(x)$ decreases.

Furthermore, since $e^{-s}e^{x(1-e^{-s/x})}\geq e^{-s}$, we get from
Lemma \ref{lemma:key}, 
$$
2\int_{0}^{x\varphi (x)}e^{-s}\, ds \leq \frac{2}{3}+
\int_{1}^{\infty}\xi'(t)\, dt = \frac{2}{3}+ \lim _{x\to \infty}\xi(x)-\xi(1)=1,
$$
since $\lim _{x\to \infty}\xi(x)=1$. (See
Proposition \ref{prop:xi_properties}.) From this we obtain that
$x\varphi(x)$ is bounded by $\log 2$ for all $x>0$. 

If we let $l=\lim_{x\to 0_+}x\varphi(x)$ then by the dominated
convergence theorem we find
$$
2\int_{0}^{l}e^{-s}\, ds = \frac{2}{3}+
\int_{1}^{\infty}\xi'(t)\, dt = 1,
$$
so that $l=\log 2$.

We let finally $L=\lim_{x\to \infty}x\varphi(x)$. We get in the same
way as before
$$
2\int_{0}^{L}1\, ds = \frac{2}{3},
$$
so that $L=1/3$. \hfill $\square$

{\it Proof of Theorem \ref{thm:cr_phi}.} 
We use the expressions $A(x)$ and $B(x)$ and the relation
(\ref{eq:smart}) from the proof of Proposition \ref{prop:xphi}. We get
\begin{eqnarray*}
1-x\varphi'(x) & = &
e^{x(\varphi(x)-1+e^{-\varphi(x)})}(A(x)+B(x))+\varphi(x)+1\\
&\leq & e^{x\varphi(x)}(A(x)+B(x))+\varphi(x)+1.
\end{eqnarray*}
Now, using that $1-e^{-a}<a$ and $1-(1+a)e^{-a}< a^2/2$  for
$a> 0$, we find
$$
A(x)\leq \int_{0}^{x\varphi (x)}e^{-s}e^{x(s/x)}\left(
  \frac{s^2}{2x^2}\right)\, ds =
\frac{x\varphi(x)^3}{6}.
$$
Furthermore, since $\xi'(t+1)\leq \xi'(1) = 8/135$ for $t\geq 0$ (see again
Proposition \ref{prop:xi_properties}),
$$
B(x)\leq \frac{1}{2}\frac{8}{135}\int_{0}^{\infty}te^{-xt}\, dt =
\frac{4}{135x^2}.
$$
This gives
$$
1-x\varphi'(x)\leq 1+\varphi(x)+e^{x\varphi(x)}\left(
  \frac{x\varphi(x)}{6}+\frac{4}{135(x\varphi(x))^2}\right)
\varphi(x)^2.
$$
It is easily seen that the function 
$$
\rho (u)=e^{u}\left(
  \frac{u}{6}+\frac{4}{135u^2}\right)
$$
attains its maximum on $[1/3, \log 2]$ at $u=1/3$ with value 
$$
\rho \left(\frac{1}{3}\right)=\frac{29}{90}e^{1/3}<\frac{1}{2}.
$$
We obtain from this the relation
$$
1-x\varphi'(x)<1+\varphi(x)+\frac{1}{2}\varphi(x)^2 <e^{\varphi(x)}.
$$
\hfill $\square$

\begin{rem}
The difference between $e^{\varphi(x)}$ and
$1+\varphi(x)+\frac{1}{2}\varphi(x)^2$ is $O(\varphi(x)^3)$, that is
(by Proposition \ref{prop:xphi}) the difference is $O(x^{-3})$. Hence
for large $x$ the difference is very small, reflecting the fact that
$m'(x)$ is very close to 1. For $x$ close to 0, the difference is
large, reflecting the fact that $m'(x)$ approaches 0 rapidly. See
Proposition \ref{prop:m_at_zero}.
\end{rem}

\section{Asymptotic results}
\label{sec:asymptotic}
Here we describe the behaviour of $m$ and $m'$ at zero and at
infinity. We first investigate $m$ near zero. We use $f(x)\sim g(x)$
for $x\to 0$ to denote that 
$$
\lim_{x\to 0}f(x)/g(x) = 1.
$$
\begin{prop}
\label{prop:m_at_zero}
We have $m(x)^x\to 1/2$ as $x \to 0$ and 
$$
m(x)\sim e^{-\gamma}2^{-1/x} \quad \text{as} \quad
x\to 0.
$$
Furthermore,
$$
m'(x)\sim (\log 2)e^{-\gamma}\frac{2^{-1/x}}{x^2} \quad \text{as} \quad
x\to 0.
$$
\end{prop}
{\it Proof.} For the
function $l(x)\equiv \log (m(x)^x)$ we have $l(x)=-x\varphi (x)+x\log
x$. Since $x\log x\to 0$ as $x\to 0$ then $l(x)\to -\log 2$ as
$x\to 0$. Therefore $m(x)^x\to 1/2$ as $x \to 0$.

By definition of $m(x)$ and by the functional
equation of the gamma function we have
$$
\int_0^{m(x)}e^{-t}xt^{x-1}\, dt=\frac{1}{2} \Gamma(x+1).
$$
We perform integration by parts on the integral on the left hand side
in this relation and we thus get
$$
m(x)^xe^{-m(x)}+\int_0^{m(x)}e^{-t}t^{x}\, dt=\frac{1}{2} \Gamma(x+1).
$$
We next differentiate this relation and get after some manipulation
$$
\log m(x) +\frac{x}{m(x)}\, m'(x) = e^{m(x)}m(x)^{-x}\left( \frac{1}{2}
  \Gamma'(x+1) -\int_0^{m(x)}(\log t)\, e^{-t}t^{x}\, dt\right) .
$$
We now use that $m(x)^x\to 1/2$ as $ x\to 0$ and we get in this way
\begin{equation}
\label{eq:gamma_relation}
\lim_{x\to 0}\left( \log m(x) +\frac{x}{m(x)}\, m'(x)\right) = \Gamma
'(1)= -\gamma.
\end{equation}
This is the same as $l'(x)=(\log (m(x)^x))'\to -\gamma$ as $x\to
0$. Using l'Hospital's rule we get
$$
\frac{l(x)+\log 2}{x}\rightarrow -\gamma
$$
for $x\to 0$. Therefore $m(x)\sim 2^{-1/x} e^{-\gamma}$ as $x\to
0$.

By (\ref{eq:gamma_relation}) we get 
$$
x\log m(x)+\frac{x^2}{m(x)}m'(x)\to 0
$$
for $x\to 0$ and, since (as used before) $x\log m(x)\to -\log 2$, we
get
$$
m'(x)\sim \frac{\log 2}{x^2}\, 2^{-1/x} e^{-\gamma}
$$
as $x\to 0$. \hfill $\square$

At infinity we have
\begin{prop}
\label{prop:m_at_infinity}
The functions $m$ and $m'$ have asymptotic expansions at infinity. We
have in particular 
$$
m'(x)=1-\frac{8}{405x^2}-\frac{368}{25515x^3}+o(x^{-3}) \quad \text{as} \quad
x\to \infty
$$
and 
$$
m(x)=x-\frac{1}{3}+\frac{8}{405x}+\frac{184}{25515x^2}+o(x^{-2}) \quad \text{as} \quad
x\to \infty.
$$
\end{prop}

\begin{rem}
  Choi (\cite{choi}) found the asymptotic expansion of $m(n+1)$ up to
  order $o(n^{-3})$. Higher order expansions of $m(n+1)$ were found in
  \cite{good}. In the appendix we have included higher order
  expansions of $m(x)$ and $\varphi(x)$ (and higher order derivatives
  of the function $\xi$ at 1). Because of the complexity, the computations behind
  the expansions in the appendix were made using ``Maple 9'', using
  the same method as described in Lemma \ref{lemma:phi_asymptotic}.
\end{rem} 
 
\begin{lemma}
\label{lemma:phi_asymptotic}
The functions $\varphi$ and $\varphi'$ have asymptotic expansions at
infinity. For the function $\varphi$ we have in particular
$$
\varphi(x)=\frac{1}{3x}+\frac{29}{810x^2}-\frac{37}{25515x^3}+o(x^{-3})
$$
as $x\to \infty$.
\end{lemma}
{\it Proof.} We have already seen in Proposition \ref{prop:xphi} that
$\varphi(x)=O(x^{-1})$. Therefore 
$$
2\int_0^{\varphi(x)}e^{-x(u+e^{-u}-1)}\, du =
2\sum_{k=0}^{n}\frac{(-1)^kx^k}{k!}\int_0^{\varphi(x)}(u+e^{-u}-1)^k\, du +
o(x^{-(n+1)}).
$$ 
On the other hand, by partial integration,
\begin{equation}
\label{eq:partial_integration}
\int_{0}^{\infty}\xi(t+1)e^{-xt}\, dt=
\sum_{k=0}^{n}\frac{\xi^{(k)}(1)}{x^{k+1}} + o(x^{-(n+1)}).
\end{equation}
(We remark here that the sum involving the derivatives of $\xi$ at 1
describes the asymptotic behaviour of Ramanujans function, see Section
\ref{sec:ramanujan}.)
From Proposition \ref{prop:key} we thus get
\begin{equation}
\label{eq:phi_asymptotic}
\varphi(x)= 
\sum_{k=0}^{n}\frac{\xi^{(k)}(1)}{2x^{k+1}}- \sum_{k=1}^{n}\frac{(-1)^kx^k}{k!}\int_0^{\varphi(x)}(u+e^{-u}-1)^k\, du + o(x^{-(n+1)}).
\end{equation}
From this relation it is possible to deduce that $\varphi(x)$ has an asymptotic
expansion of the form 
\begin{equation}
\label{eq:phi_asymptotic_n}
\varphi(x)= 
\sum_{k=1}^{n}\frac{c_k}{x^k}+o(x^{-n})
\end{equation}
as $x\to \infty$ and for any $n$. (See below.)

Let us first find the coefficients $c_1,c_2$
and $c_3$ in the expansion. If $n=0$ we get, using $\xi(1)=2/3$,
$$
\varphi(x)= \frac{1}{3x}+ o(x^{-1}).
$$
(This we have already found in Proposition \ref{prop:xphi}.) For $n=1$ we get, using $\xi'(1)=8/135$, 
$$
\varphi(x)= \frac{1}{3x}+\frac{4}{135x^{2}}+
x\int_0^{\varphi(x)}(u+e^{-u}-1)\, du + o(x^{-2}).
$$ 
Here we use $u+e^{-u}-1= u^2/2+o(u^2)$ for $u\to 0$ and we get in this
way
$$
\varphi(x)= \frac{1}{3x}+\frac{4}{135x^{2}}+
x\left( \frac{\varphi(x)^3}{6}+o(x^{-3})\right) +o(x^{-2}).
$$
Then we use the result for $n=0$ to obtain that
$x\varphi(x)^3=1/(27x^2)+o(x^{-2})$, and if we insert this into the
relation above we get
$$
\varphi(x)= \frac{1}{3x}+\left( \frac{4}{135}+\frac{1}{162}\right) \frac{1}{x^2}+o(x^{-2})= \frac{1}{3x}+\frac{29}{810x^2}+o(x^{-2}).
$$
We repeat the argument to obtain the following for $n=2$:
\begin{eqnarray*}
\varphi(x) & = &\frac{1}{3x}+\frac{4}{135x^{2}}-\frac{8}{2835x^{3}}\\
&& +x\int_0^{\varphi(x)}(u+e^{-u}-1)\, du -\frac{x^2}{2}\int_0^{\varphi(x)}(u+e^{-u}-1)^2\, du + o(x^{-3}),
\end{eqnarray*} 
where we use $u+e^{-u}-1= u^2/2-u^3/6+o(u^3)$ and $(u+e^{-u}-1)^2=
u^4/4+o(u^4)$ for $u\to 0$ and we get
\begin{eqnarray*}
\varphi(x) & = & \frac{1}{3x}+\frac{4}{135x^{2}}-\frac{8}{2835x^{3}}\\
&&+x\left(
  \frac{\varphi(x)^3}{6}-\frac{\varphi(x)^4}{24}+o(x^{-4})\right)\\
&&-\frac{x^2}{2}\left( \frac{\varphi(x)^5}{20}+o(x^{-5})\right) + o(x^{-3}).
\end{eqnarray*} 
In the term $x\varphi(x)^3/6$
we substitute the expansion of $\varphi(x)$ for $n=1$ and in the terms
$x \varphi(x)^4/24$ and $x^2\varphi(x)^5/40$
the expansion of $\varphi(x)$ for $n=0$. In this way we get
\begin{eqnarray*}
\varphi(x) & = & \frac{1}{3x}+\frac{29}{810x^2}+\left(
  -\frac{8}{2835}+\frac{1}{18}\frac{29}{810}-\frac{1}{24}\frac{1}{81}-\frac{1}{40}\frac{1}{243}\right) \frac{1}{x^3} + o(x^{-3})\\
& = & \frac{1}{3x}+\frac{29}{810x^2}-\frac{37}{25515x^3}+o(x^{-3}).
\end{eqnarray*}
These computations also indicate how to show that there is an
asymptotic expansion (\ref{eq:phi_asymptotic_n}) for every $n\geq
1$. One could use an inductive argument based on the relation
(\ref{eq:phi_asymptotic}): First of all, the sum 
$$
\sum_{k=0}^{n}\frac{\xi^{(k)}(1)}{2x^{k+1}}
$$
contains a term of order $1/x^{n+1}$. Next the integrand $(u+e^{-u}-1)^k$ is
approximated by its Taylor polynomial of order $n+k$,
$$
(u+e^{-u}-1)^k=\sum_{l=2k}^{n+k}\alpha_{k,l}u^l+ o(u^{n+k}),
$$
and the expansion
$\varphi(x)=\sum_{k=1}^nc_k/x^{-k}+o(x^{-n})$ is then used in the
upper limit of the integrals in the sum
$$
\sum_{k=1}^{n}\frac{(-1)^kx^k}{k!}\int_0^{\varphi(x)}(u+e^{-u}-1)^k\,
du.
$$
Using these approximations it is possible to obtain
$$
\varphi(x)= 
\sum_{k=0}^{n}\frac{\xi^{(k)}(1)}{2x^{k+1}}-
\sum_{k=1}^{n}\frac{(-1)^kx^k}{k!}\sum_{l=2k}^{n+k}\frac{\alpha_{k,l}}{l+1}\left(
  \sum_{j=1}^n\frac{c_j}{x^j}\right)^{l+1} + o(x^{-(n+1)}).
$$
This is an expansion of $\varphi(x)$ of order $n+1$.

To see that also $\varphi'(x)$ has an asymptotic expansion we
differentiate (\ref{eq:key}) and get
$$
2\varphi'(x)e^{-x(\varphi(x)+e^{-\varphi(x)})}-2\int_0^{\varphi(x)}e^{-x(u+e^{-u})}(u+e^{-u})\,
  du =-\int_1^{\infty}t\xi(t)e^{-xt}\, dt.
$$
Adding (\ref{eq:key}) to the relation above we get after multiplication by $e^x$ and a
change of variable $s=t-1$,
\begin{eqnarray*}
\lefteqn{2\varphi'(x)e^{-x(\varphi(x)+e^{-\varphi(x)}-1)}}\\
& = & 2\int_0^{\varphi(x)}e^{-x(u+e^{-u}-1)}(u+e^{-u}-1)\,
  du-\int_0^{\infty}s\xi(1+s)e^{-xs}\, ds.
\end{eqnarray*}
In the first integral we consider again Taylor approximation and in
the second integral we perform integration by parts. Using also the
asymptotic expansion of $e^{x(\varphi(x)+e^{-\varphi(x)}-1)}$ we are
finally able to see that there is an asymptotic expansion of
$\varphi'(x)$. The coefficients in this expansion can be identified by
integrating the expansion and using the known expansion of
$\varphi(x)$, cf \cite{A:A:R}[Appendix C].
\hfill $\square$

{\it Proof of Proposition \ref{prop:m_at_infinity}.} Since
$\varphi(x)$ has an asymptotic expansion,
$$
m(x)=xe^{-\varphi(x)}=x\sum_{k=0}^{\infty}\frac{(-1)^k\varphi(x)^k}{k!}
$$ 
also has an asymptotic expansion.  We have in particular 
$$
m(x)=xe^{-\varphi(x)}=x\left(1-\varphi(x)+\frac{\varphi(x)^2}{2}-\frac{\varphi(x)^3}{6}+o(x^{-3})\right).
$$
We insert in this relation the asymptotic expansion of $\varphi(x)$
from the lemma above. We get, after some computation,
\begin{eqnarray*}
m(x) & = & x\left(1-\frac{1}{3x}+\frac{8}{405x^2}+\frac{184}{25515x^3}+o(x^{-3})\right)\\
& = & x-\frac{1}{3}+\frac{8}{405x}+\frac{184}{25515x^2}+o(x^{-2}).
\end{eqnarray*}
Since $\varphi'(x)$ has an asymptotic expansion, the same is true for
$m'(x)$, since $m'(x)=(1-x\varphi'(x))e^{-\varphi(x)}$. The expansion
of $m(x)$ can be found by integrating the expansion of $m'(x)$, and
this gives the desired expansion of $m'(x)$. \hfill $\square$.

\begin{rem}
It is clear that the methods above can be continued so the asymptotic
behaviour of $m^{(k)}(x)$, $k\geq 2$ for $x\to 0$ and $x\to \infty$
can be determined by differentiation of the asymptotic formulas for
$m'(x)$.
\end{rem}

\section{Properties of the auxiliary function $\xi$}
\label{sec:xi}
In this section we derive the properties of the function $\xi$,
stated in Proposition \ref{prop:xi_properties}.
Most of these properties we found surprisingly difficult to establish. Much of
the difficulty lies in the fact that $\xi$ in given as a sum of two
terms, where it is necessary to control the cancellation between these
two terms. 

Throughout this section $u$ and $v$ denote the functions defined in
(\ref{eq:xi}). We begin by considering them separately. The investigation is based on a simple lemma
of independent interest. We recall that a $C^{\infty}$-function $\nu
(t)$ defined on the positive half-line is called {\em completely monotonic} if 
$$
(-1)^p\nu ^{(p)}(t)\geq 0, \quad \text{for all integers} \quad p\geq 0.
$$
It is called a {\em Bernstein function} if
$$
\nu(t)\geq 0 \quad \text{and} \quad(-1)^p\nu ^{(p)}(t)\leq 0, \quad \text{for all integers} \quad p\geq 1.
$$
The last conditions can also be expressed that $\nu'(t)$ is completely
monotonic. For details about these classes of functions see e.g.\ \cite{berg-forst}.
\begin{lemma}
Let $F$ be completely monotonic, let $\varsigma(t)$ be a positive
$C^{\infty}$-function for $t>0$ and assume that $\varsigma'(t)=F(\varsigma(t))$.
Then we have:

For each $n\geq 1$ there exists a completely monotonic function $F_n$
such that
$$
(-1)^{n-1}\varsigma ^{(n)}(t)=F_n(\varsigma(t)).
$$
In particular $\varsigma$ is a Bernstein function.
\end{lemma}
{\it Proof.} For $n=1$ we can use $F_1=F$. 
If $(-1)^{n-1}\varsigma ^{(n)}(t)=F_n(\varsigma(t))$ 
for some completely monotonic function $F_n$ then
$$
(-1)^{n}\varsigma ^{(n+1)}(t)=-F_n'(\varsigma (t))\varsigma '(t)=
-F_n'(\varsigma (t))F(\varsigma (t)).
$$
Here $F_{n+1}\equiv -F_n'F$ is completely monotonic as a product of
two completely monotonic fuctions. 

Furthermore we have
$$
(-1)^{n}\varsigma ^{(n)}(t)= -F_n(\varsigma(t)) \leq 0,
$$
so $\varsigma$ is a Bernstein function. \hfill $\square$ 

\begin{prop}
\label{prop:u_and_v}
The functions $-u(t+1)$ and $v(t+1)$ are Bernstein functions and 
$$
\lim_{t\to
  \infty}v(t)/t=1, \quad \lim_{t\to \infty}u(t)/t=0.
$$
\end{prop}
{\it Proof.} We have $v'(t)=F(v(t))$ where 
$$
F(x)=\frac{1}{1-e^{-x}}=\sum_{k=0}^{\infty}e^{-kx}
$$
is completely monotonic. Since $v(t)\to \infty$ as $t\to \infty$ we
have $e^{-v(t)}\to 0$ and hence $v(t)/t=1-e^{-v(t)}/t \to 1$ as $t\to
\infty$.

For the function $w\equiv -u$ we have $w'(t)=G(w(t))$ where
$$
G(x)=\frac{1}{e^{x}-1}=\sum_{k=1}^{\infty}e^{-kx} 
$$ 
is completely monotonic. We have $t=e^{-u(t)}+u(t)>1+u(t)^2/2$ (since
$u(t)<0$) so that $-u(t)<\sqrt{2(t-1)}$ and thus $u(t)/t\to 0$ as
$t\to \infty$. \hfill $\square$

Bernstein functions admit integral representations (see e.g.\
\cite[p.\ 64]{berg-forst}) and thus we have
$$
v(t+1)=t+\int_0^{\infty}(1-e^{-xt})\, d\lambda (x)
$$
and 
$$
-u(t+1)=\int_0^{\infty}(1-e^{-xt})\, d\sigma (x),
$$
for some positive measures $\lambda$ and $\sigma$ on $(0,\infty)$.
Since $t+1-v(t+1)=e^{-v(t+1)}\to 0$ for $t \to \infty$ we conclude that
$\lambda ((0,\infty))=1$, hence
$$
v(t+1)=t+1-\int_0^{\infty}e^{-xt}\, d\lambda (x).
$$
The measure $\sigma$ has infinite total mass because $u(t)\to -\infty$
for $t\to \infty$.

From these representations we see that $\xi (t)=u'(t)+v'(t)\to 1$ as
$t\to \infty$ and also that $u^{(n+1)}(t)$ and
$v^{(n+1)}(t)$ both tend to zero as $t$ tends to infinity for any $n\geq
1$. This shows that $\xi^{(n)}(t)=u^{(n+1)}(t)+v^{(n+1)}(t)\to 0$ as
$t\to \infty$ for any $n\geq 1$.

One could hope to deduce further properties
of $\xi$ by using these integral representations. We have
not succeeded in doing this, since much cancellation between $u$ and
$v$ takes place and we do not know $\lambda$ and $\sigma$ explicitly. Furthermore, it it
not even true that $\xi(t+1)$ is a Bernstein function even though it is increasing and concave.

We have found an approach using complex analysis and
we shall make extensive use of the theory of the so-called Pick
functions, see \cite{D}. A holomorphic function $p$ defined in the upper
half plane is a Pick function if it maps the upper half plane into
the closed upper half plane, or put in another way, if $\Im p$ is a non-negative
harmonic function. A Pick function has an integral representation, 
\begin{equation}
\label{eq:general_pick}
p(z)=az+b+\int_{-\infty}^{\infty}\left(
  \frac{1}{t-z}-\frac{t}{t^2+1}\right) \, d\mu (t),
\end{equation}
where $a\geq 0$, $b$ is real and $\mu$ is a positive measure on
the real line. Furthermore,
$$
a=\lim_{y\to \infty}\frac{\Im p(iy)}{y}, \quad b=\Re p(i)
$$
and 
$$
\mu= \lim_{y\to 0_+}\frac{\Im p(x+iy)}{\pi}dx
$$
in the vague topology. 
Note that a Pick function can be extended to a holomorphic function in
$\mathbb C\setminus \mathbb R$ by (\ref{eq:general_pick}). If
$A\subseteq \mathbb R$ is closed  we recall that $p$ has a holomorphic
extension to $\mathbb C\setminus A$ if and only if the support of
$\mu$ is contained in $A$.

The proofs of the results to follow are based on an investigation of the
holomorphic function 
\begin{equation}
\label{eq:f}
f(z)=e^{-z}+z.
\end{equation}

Concerning the
values of the function $\xi$ and its derivatives at $z=1$ we have the following proposition. 
\begin{prop}
\label{prop:xi_at_1}
The function $\xi(z)$ is holomorphic in a neighbourhood of $z=1$ and we have $\xi(1)=2/3$ and $\xi'(1)=8/135$.
\end{prop}
{\it Proof of Proposition \ref{prop:xi_at_1}.} The function $f(z)-1$
has a zero of multiplicity 2 at $z=0$. Hence there exists a
holomorphic function $h$ in a neighbourhood of $0$ such that
$f(z)-1=h(z)^2$ there. Since $h'(0)^2=1/2\ne 0$, $h$ is one-to-one
near $z=0$. We choose $h'(0)=1/\sqrt{2}$ and with this choice $h$ is uniquely determined.

We have thus a radius $r>0$ such that for any $w$ with $|w|<r$,
there are exactly two solutions to the equation $f(z)=w$,
namely $z=h^{-1}(\pm \sqrt{w-1})$. In particular, for $t>1$ and close to
$1$ we have 
$$
u(t)=h^{-1}(-\sqrt{t-1})\quad \text{and} \quad  v(t)=h^{-1}(\sqrt{t-1}),
$$
where $u$ and $v$ are the functions appearing in (\ref{eq:xi}). We denote by $\sum_{n=1}^{\infty}a_nw^n$ the Taylor series of $h^{-1}(w)$. 
Then 
$$
u(t)+v(t)= \sum_{k=1}^{\infty}2a_{2k}(t-1)^k,
$$
and, since $\xi(t)=u'(t)+v'(t)$, we find $\xi(1)=2a_2$ and
$\xi'(1)=4a_4$. The numbers $a_2$ and $a_4$ can be found as follows. If
$h(z)=\sum_{n=1}^{\infty}b_nz^n$ then
$$
\left( \sum_{n=1}^{\infty}b_nz^n\right) ^2 = e^{-z}+z-1 = \frac{1}{2}z^2
-\frac{1}{6}z^3 +\frac{1}{24}z^4+\ldots
$$
and this gives after some manipulation $b_1=1/\sqrt{2}$, $b_2=-\sqrt{2}/12$,
$b_3=\sqrt{2}/72$ and $b_4=-\sqrt{2}/540$. The numbers $a_2$ and $a_4$
can now be identified if we differentiate four times the relation
$(h^{-1})(h(z))=z$ and put $z=0$ (so that expressions
involving e.g.\ $(h^{-1})^{(4)}(0)$ appear).
We find, again after some manipulation,
$(h^{-1})^{(2)}(0) = 2/3$ and $(h^{-1})^{(4)}(0) = 16/45$,
so that $\xi(1)=2a_2=2/3$ and $\xi'(1)=4a_4=4\cdot (16/45)\cdot
(1/24)=8/135$. \hfill $\square$

\begin{rem}
One can in principle find any derivative of $\xi$ at 1 in this
way. Computations show that $\xi''(1)=-16/2835$. The higher order
derivatives of $\xi$ at 1 can be found in the appendix.
\end{rem}

\begin{figure}[htb]
\begin{center}
\includegraphics[scale=0.6]{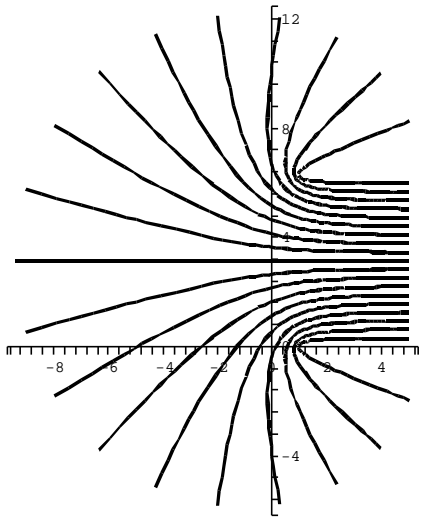}
\end{center}
\caption{The image of horizontal lines under $f$}
\label{figure:f}
\end{figure}

The monotonicity properties of $\xi$ follows from the next proposition.
\begin{prop}
\label{prop:xi_incr_conc}
The function $\xi$ increases on the real line and it is concave on
$[1, \infty)$.
\end{prop} 
The key to the proof of this proposition is an integral representation
of $\xi$ that we give in Proposition \ref{prop:xi_int}.  
\begin{lemma}
\label{lemma:f_conformal}
The function $f(z)=e^{-z}+z$ is a conformal mapping of the strip 
$$
\mathcal S_{2\pi} =\{ 0< \Im z<2\pi\}
$$ 
onto the domain 
$$
\mathcal T = \mathbb C\setminus
\left( \{ \Im w=2\pi i, \Re w\geq 1\} \cup \{ \Im w=0, \Re
  w\geq 1\} \right) .
$$
\end{lemma} 
In Figure \ref{figure:f} we have indicated the image of some
horizontal lines in $\mathcal S_{2\pi}$ under the function $f$.

{\it Proof.} We have
$f(x+iy)=\sigma +i\tau $ if and only if
\begin{eqnarray*}
e^{-x}\cos y +x & = & \sigma \\
-e^{-x}\sin y +y & = & \tau .
\end{eqnarray*}
For $y=\pi$ we get $\tau = \pi$ and $\sigma = x-e^{-x}$ so $f$ maps
the horizontal line $\mathbb R +i\pi$ onto itself and is one-to-one
there.

Note that for $0<y<\pi$ we have $\tau =-e^{-x}\sin y +y<y$ so $f(\mathcal
S)\subseteq \{ \Im w<\pi\}$.

For $\tau <\pi$ we get $e^x=(\sin y)/(y-\tau )$ and since $\tau <y$ we find
$x=\log (\sin y/(y-\tau ))$ and finally we shall look for a solution $y$
to the equation 
$$
\log \frac{\sin y}{y-\tau }+\frac{(y-\tau )\cos y}{\sin y} = \sigma .
$$
For $\tau <\pi$ we  put
\begin{equation}
\label{eq:F_tau}
F_{\tau}(y)=\log \frac{\sin y}{y-\tau }+\frac{(y-\tau )\cos y}{\sin y}.
\end{equation}
We claim the following:
\begin{enumerate}
\item For $0<\tau <\pi$, $F_{\tau}(y)$ decreases for $y\in (\tau ,\pi)$ from
  $\infty$ to $-\infty$.
\item For $\tau =0$, $F_{\tau}(y)$ decreases for $y\in (\tau ,\pi)$ from
  1 to $-\infty$.
\item For $\tau <0$, $F_{\tau}(y)$ decreases for $y\in (0,\pi)$ from
  $\infty$ to $-\infty$.
\end{enumerate}
It is easy to see that $\lim_{y\to \pi_-}F_{\tau}(y)=-\infty$, $\lim_{y\to
  \tau_+ }F_{\tau}(y)=\infty$ for $\tau >0$ and $\lim_{y\to 0_+}F_0(y)=1$. We also have
$\lim_{y\to 0_+}F_{\tau}(y)=\infty$ for $\tau <0$ since the first term in
$F_{\tau}(y)$ has a logarithmic singularity.  To verify 1., 2.\ and 3.\ we need
to show that $F_{\tau}'(y)<0$. We find, after some computation,
$$
F_{\tau}'(y)= \frac{1}{\sin y}\left( 2\cos y - \frac{\sin
    y}{y-\tau }-\frac{y-\tau }{\sin y} \right).
$$
If we put $\kappa=(\sin y)/(y-\tau )$ then $\kappa >0$ and hence $\kappa +1/\kappa
\geq 2$. Since $\cos y<1$ we see that indeed $F_{\tau}'(y)<0$.

From 1.\ and 3.\ it follows that for given $\sigma \in \mathbb R$ and
$\tau <\pi$, $\tau \neq 0$ there is a unique solution $y$ to the equation
$F_{\tau}(y)=\sigma $ and therefore there is a unique solution to $f(x+iy)=\sigma +i\tau $.
If $\tau =0$ there is by 2.\ a unique solution $y$ to the equation
$F_{\tau}(y)=\sigma $ when $\sigma <1$ and none when $\sigma \geq 1$.

For $\pi <y<2\pi$  we have $\sin y<0$ so $\tau >y$ and therefore  
$$
f(\{ \pi < \Im z <2\pi \}) \subseteq \{ \pi < \Im w  \} .
$$
For $\tau >\pi$ we put 
$$
F_{\tau}(y)=\log \frac{-\sin y}{\tau -y}+\frac{(\tau -y)\cos y}{-\sin y}.
$$
It follows that $F_{\tau}(y)$ increases with $y\in (\pi ,2\pi )$, and
the conformality follows in the same way as for the case $\tau <\pi$.
\hfill $\square$

\begin{lemma}
\label{lemma:f_pick}
The function $w\mapsto f^{-1}(1-w)$ is a Pick function with the
representation 
$$
f^{-1}(1-w)= \Re f^{-1}(1-i) + \int_0^{\infty}\left(
  \frac{1}{t-w}-\frac{t}{t^2+1}\right) \eta (t)\, dt,
$$
where  
$$
\eta (t)=\frac{1}{\pi}\, \Im f^{-1}(1-t)
$$ 
is increasing on $(0, \infty)$ from 0 to 1. 
\end{lemma}
{\it Proof.} First of all, $1-w$ belongs to the lower half plane when
$w$ belongs to the upper half plane. Since $f^{-1}$ maps all of
$\mathcal T$, and hence in particular the lower half plane, into the
strip $0<\Im z <2\pi$, $f^{-1}(1-w)$ is certainly a Pick function. We next derive
its integral representation. Since its imaginary part is bounded the
number $a$ in the representation (\ref{eq:general_pick}) must be zero. 

It remains to identify the measure $\mu$. Since the function
$f^{-1}(1-t-iy)$ is continuous on e.g.\ $\mathbb R \times [0,1]$ and
is real for $y=0$ and $t<0$ we find 
$$
\eta (t)=\lim_{y\to 0_+}\frac{\Im f^{-1}(1-t-iy)}{\pi}=\frac{\Im f^{-1}(1-t)}{\pi}
$$
for $t>0$ and $\eta(t)=0$ for $t<0$.

By definition the function $Y(t)=\pi \eta(t)$ satisfies the equation
$$
F_0(Y(t))=\log \frac{\sin Y(t)}{Y(t)}+\frac{Y(t)\cos Y(t)}{\sin
  Y(t)}=1-t
$$ 
(with the notation of Lemma \ref{lemma:f_conformal}). We find from
this
$$
Y'(t)F_0'(Y(t))=-1,
$$ 
and since $F_0'$ is negative (see again Lemma
\ref{lemma:f_conformal}), $Y(t)$ and hence $\eta(t)$ must be
increasing. Since $Y(t)$ tends to $\pi$ as $t$ tends to $\infty$ (this
is because $F_0(Y(t))\to -\infty$),
$\eta(t)\to 1$ as $t\to \infty$.\hfill $\square$

\begin{lemma}
\label{lemma:g_pick}
The function 
$$
g(w)= f^{-1}(w)+ \overline{f^{-1}(\overline{w})}
$$
has a holomorphic extension to $\mathbb C\setminus
 \{ \Im w=\pm 2\pi i, \Re w\geq 1\}$. 

For $t>1$ we have $g(t)=u(t)+v(t)$ and $g'(t)=\xi(t)$. The function $g'$ is thus a holomorphic extension of
$\xi$.

All points on the lines $\{ \Im
 w=\pm 2\pi i, \Re w\geq 1\}$ are singular points for $g$.
\end{lemma}

{\it Proof.} Since $f^{-1}$ is a conformal mapping of the region 
$$
\mathbb C\setminus
 \left( \{ \Im w=2\pi i, \Re w\geq 1\} \cup \{ \Im w=0, \Re w\geq
   1\}\right) ,
$$
the function $g$ is holomorphic in 
$$
\mathbb C\setminus
 \left( \{ \Im w=\pm 2\pi i, \Re w\geq 1\} \cup \{ \Im w=0, \Re w\geq
   1\}\right).
$$
However, $g$ has a continuous extension to $\{ \Im w=0, \Re w\geq
   1\}$ with boundary values $u+v$, so from Moreras theorem we
   conclude that $g$ is holomorphic across this half-line.  Concerning
   the boundary values on the half-lines $\{ \Im w=\pm 2\pi i, \Re
   w\geq 1\}$ we have ($\sigma \geq 1$)
\begin{eqnarray*}
\lim_{\tau \to 2\pi _-}g(\sigma+i\tau) & = &  \lim_{\tau \to 2\pi
  _-}f^{-1}(\sigma+i\tau) +\overline{f^{-1}(\sigma-2\pi i)}\\
& = &  v(\sigma)+2\pi i +\overline{f^{-1}(\sigma-2\pi i)}
\end{eqnarray*}
and similarly
$$
\lim_{\tau \to 2\pi _+}g(\sigma+i\tau)=  u(\sigma)+2\pi i
+\overline{f^{-1}(\sigma-2\pi i)}. 
$$
Therefore $g$ is not continuous across any segment of the half-line 
$$
\{ \Im w=2\pi i, \Re w\geq 1\}
$$
and so all these points are singular points. Since
   $g(w)=\overline{g(\overline{w})}$ the same conclusion holds for the
   points on the other half-line.
\hfill $\square$ 
\begin{prop}
\label{prop:divergence}
The Taylor series for $\xi(z)$ centered at 1
has radius of convergence equal to $2\pi$ and the asymptotic series
$$
\sum_{k=0}^{\infty}\frac{\xi^{(k)}(1)}{z^k}
$$
diverges for any $z$ in $\mathbb C$.
\end{prop}
{\it Proof.} We have $\xi = g'$, where $g$ is the function in Lemma
\ref{lemma:g_pick}. If the radius of convergence of $\xi$ at 1 were
larger than $2\pi$ then the primitive $g$ would also have a holomorhic
extension to this larger disk, and this contradicts the fact that $g$
has singular points in that disk.

Concerning the divergence of the asymptotic series we use that the
radius of convergence of the Taylor series is finite. It means that
$$
\limsup_{k\to \infty}\left( \frac{|\xi^{(k)}(1)|}{k!}\right) ^{1/k}
>0,
$$ 
and hence that for some $\epsilon>0$,
$|\xi^{(k)}(1)|^{1/k}/(k!)^{1/k}>\epsilon$ for infinitely many
$k$. From Stirlings formula we have $(k!)^{1/k} \sim k/e$ and
therefore
$$
\limsup_{k\to \infty}|\xi^{(k)}(1)|^{1/k} =
\infty.
$$
This shows on the other hand that the asymptotic series
diverges for any complex number $z$. \hfill $\square$

In the following $\log$ denotes the principal logarithm defined in the cut
plane $\mathbb C\setminus (-\infty, 0]$.
\begin{lemma}
\label{lemma:psi_pick}
The function 
$$
\Psi(w) = g(\log w) = f^{-1}(\log w)+ \overline{f^{-1}(\log \overline{w})}
$$
is a Pick function. It has the representation
$$
\Psi(w)=\Re
\Psi(i)-\int_{0}^{\infty}\left(\frac{1}{t+w}-\frac{t}{t^2+1}\right) \,
h(t)dt,
$$
where the function $h$ is given as 
$$
h(t)=1-\frac{\Im \, f^{-1}(\log t-i\pi)}{\pi}.
$$
Furthermore, $h$ increases on $(0,\infty)$ from 0 to 1.

For $\tau \in [e,\infty)$ we have 
$$
\Psi(\tau )=u(\log \tau)+v(\log \tau).
$$
\end{lemma}

{\it Proof.} We consider the holomorphic function $g$ from Lemma
\ref{lemma:g_pick} in the strip $\mathcal S=\{ 0<\Im w<\pi\}$.
Hence the function $V(w)=\Im g(w)$ is a harmonic function in $\mathcal S$. Furthermore,
since $f^{-1}$ maps all of $\{ \Im w<\pi\}\setminus \{ x\, | \, x\geq 1\}$ into $\mathcal S$, $V$ is also
bounded there, with the apriori bound $-\pi <V(w) <\pi$. We
claim that indeed $0<V(w)<\pi$ for all $w\in \mathcal S$.

We consider the boundary values of $V$. The horizontal line
$\{t+i\pi, \, |\, t \in \mathbb R\}$ is mapped by $f^{-1}$ to itself,
whereas $\{t-i\pi, \, |\, t \in \mathbb R\}$ is mapped by $f^{-1}$ to
some curve inside the strip $\mathcal S$. Therefore $\overline{f^{-1}(t-i\pi)}$
has imaginary part greater than $-\pi$ and hence the boundary values
$V(t+i\pi), t\in \mathbb R$, are all non negative. Since $g(\overline{w})=\overline{g(w)}$,
the boundary values $V(t), t\in \mathbb R$, are all zero. We conclude that $V$ has
non-negative boundary values. Since it is bounded, the maximum
principle in an unbounded region (see e.g.\ \cite{Koosis}) yields that $V(w)>0$ for all $w\in
\mathcal S$.

Since $\log$ maps the upper half plane onto the strip $\mathcal S$,
the function $\Psi(w)=g(\log w)$ maps the upper half plane into
itself, and is hence a Pick function. Since $\Im \Psi$ is
bounded, $\Psi$ has an integral representation of the
form
$$
\Psi(w)=\Re \Psi (i)+\int_{-\infty}^{\infty}\left(\frac{1}{t-w}-\frac{t}{t^2+1}\right) \,
d\mu (t),
$$
for some positive measure $\mu$. Here $\Im \Psi(t)=0$ for $t>0$ since
$V$ is zero on the real line. We thus find the measure $\mu$ to be
supported on the negative half-line with density
$$
1-\frac{\Im \, f^{-1}(\log (-t)-i\pi)}{\pi}, \quad t<0.
$$
After making a change of variable ($t\mapsto -t$) in the integral we conclude that
\begin{equation}
\label{eq:psi_int}
\Psi(w)=\Re
\Psi(i)-\int_{0}^{\infty}\left(\frac{1}{t+w}-\frac{t}{t^2+1}\right) \,
h(t)dt,
\end{equation}
where $h$ is the function in the statement of the lemma.

By definition of $h$, $\pi (1-h(t))= \Im f^{-1}(\log t -i\pi)$, so $h$
is increasing if the solution $Y(t)$ to the equation
$F_{-\pi}(Y(t))=\log t$ is decreasing. This is indeed the case, since
$$
F_{-\pi}'(Y(t))Y'(t)=\frac{1}{t}>0,
$$
and $F_{-\pi}'<0$.\hfill $\square$

\begin{rem}
Since it is easily verified that $V(w)=\Im g(w)$ is positive for $\pi
<\Im w<2\pi$ it follows that $g(2\log w)$ is also a Pick function.
\end{rem}

\begin{prop}
\label{prop:xi_int}
The function $\xi(\log w)$ is holomorphic in the
cut plane $\mathbb C \setminus (-\infty, 0]$ with the representation 
$$
\xi(\log w)= 1-\int_0^{\infty}\frac{t}{t+w}h'(t)\, dt,
$$
where $h(t)$ is given in Lemma \ref{lemma:psi_pick}. 
\end{prop}
If we let $w\to 0_+$ we see that
$$
\xi(\log w)\to 1-\int_0^{\infty}h'(t)\, dt = 1-\lim_{t\to
  \infty}h(t)+h(0)= 0.
$$
Hence $\xi(t)\to 0$ as $t\to -\infty$.

{\it Proof of Proposition \ref{prop:xi_int}.} As we have noticed above, $\Psi(w)=u(\log w)+v(\log w)$, for $w\geq
e$ , where $\Psi$ is the function
in Lemma \ref{lemma:psi_pick}. Hence $\xi(\log w)= w\Psi'(w)$ and so we get from (\ref{eq:psi_int}),
$$
\xi(\log w)= \int_0^{\infty}\frac{wh(t)}{(t+w)^2}\, dt.   
$$  
If we perform integration by parts on the right-hand side of this
relation then we get
$$
\xi(\log w)= \int_0^{\infty}\frac{w}{t+w}\, h'(t)\, dt,  
$$  
or
$$
\xi(\log w)= 1-\int_0^{\infty}\frac{t}{t+w}\, h'(t)\, dt.  
$$  
\hfill $\square$

{\it Proof of Proposition \ref{prop:xi_incr_conc}.} We get by differentiation of
the formula in Proposition \ref{prop:xi_int} that
$$
\frac{\xi'(\log w)}{w}= \int_0^{\infty}\frac{t}{(t+w)^2}\, h'(t)\, dt,  
$$  
and since we know that $h'(t)>0$ for all $t\in \mathbb R$, we see that
$\xi'(t)>0$ for all $t\in \mathbb R$.  

It is more technical to establish that $\xi(t)$ is concave for $t\geq 1$.
We differentiate both sides of the relation above once more and we get in this way
$$
\xi''(\log w)-\xi'(\log w)=-2w^2 \int_0^{\infty}\frac{t}{(t+w)^3}\,
h'(t)\, dt,
$$
so that, again using the integral representation of $\xi'(\log w)$,
\begin{eqnarray*}
\xi''(\log w) & = & \int_0^{\infty}\left(
  \frac{-2w^2}{(t+w)^3}+\frac{w}{(t+w)^2} \right)
th'(t)\, dt\\
& = & \int_0^{\infty} \frac{t-w}{(t+w)^3} \, wth'(t)\, dt\\
& = & w\int_0^{\infty} \frac{s-1}{(s+1)^3} \, sh'(sw)\, ds.
\end{eqnarray*}
This last integral we split into two, one for $s<1$ and another for
$s>1$. We thus get, after making the substitution $s\to 1/s$ in the
latter one,
\begin{eqnarray*}
\frac{\xi''(\log w)}{w} & = & \int_0^{1} \frac{s-1}{(s+1)^3} \,
sh'(sw)\, ds +\int_0^{1} \frac{1/s-1}{(1/s+1)^3} \,
\frac{1}{s}h'\left( \frac{w}{s} \right) \, \frac{ds}{s^2}\\ 
& = & \int_0^{1} \frac{1-s}{(1+s)^3} \,
\left( \frac{1}{s}h'\left( \frac{w}{s} \right) -sh'(sw) \right) \, ds.
\end{eqnarray*}
From this relation we see that $\xi$ is concave on $[\log w_0,\infty)$
if
\begin{equation}
\label{eq:h_inequality}
\frac{1}{s}h'\left( \frac{w}{s} \right) -sh'(sw) \leq 0 \quad \text{for}
\quad s\in (0,1) \quad \text{and} \quad w\geq w_0.
\end{equation}
Now, as noted in the proof of Lemma \ref{lemma:psi_pick},
$h(s)=1-Y_{-\pi}(\log s)/\pi$, where $Y_{-\pi}(t)\in (0,\pi)$ is the solution to
the equation $F_{-\pi}(Y_{-\pi}(t))=t$, see (\ref{eq:F_tau}). Therefore $h'(s)=-Y_{-\pi}'(\log
s)/(\pi s)$ and hence
$$
\frac{1}{s}h'\left( \frac{w}{s} \right)= -\frac{Y_{-\pi}'(\log w -\log
  s)}{\pi w}
$$
and
$$
sh'(sw) = -\frac{Y_{-\pi}'(\log w +\log
  s)}{\pi w}.
$$
We see that (\ref{eq:h_inequality}) holds provided
$$
Y_{-\pi}'(\log w -\log s) \geq Y_{-\pi}'(\log w +\log s).
$$
Since $\log s$ runs through $(-\infty, 0)$ when $s\in (0,1)$, this
condition is the same as 
$$
Y_{-\pi}'(\log w +t) \geq Y_{-\pi}'(\log w -t)
$$
for all $t\geq 0$ and all $w\geq w_0$. We aim at proving this for
$w_0=e$, so what we really should prove is the following:
\begin{equation}
\label{eq:y_inequality}
Y_{-\pi}'(\alpha +t) \geq Y_{-\pi}'(\alpha -t)
\end{equation}
for all $t\geq 0$ and all $\alpha \geq 1$. 
We verify this inequality by using Lemma \ref{lemma:f_pick}. According
to that lemma, with $w=1-\alpha -t+i\pi$, we
  have
\begin{eqnarray*}
Y_{-\pi}(\alpha +t) & = & \Im f^{-1}(\alpha +t-i\pi))\\
& = & \Im f^{-1}(1-(1-\alpha -t+i\pi))\\
& = & \int_0^{\infty}\frac{\pi}{(s+t+\alpha -1)^2+\pi^2}\, \eta(s) \, ds
\end{eqnarray*}
so that
\begin{eqnarray*}
Y_{-\pi}'(\alpha +t) 
&=&
\int_0^{\infty}\left( \frac{\partial}{\partial
    t}\frac{\pi}{(s+t+\alpha -1)^2+\pi^2}\right) \, \eta(s) \, ds\\
&=&
\int_0^{\infty}\left( \frac{\partial}{\partial
    s}\frac{\pi}{(s+t+\alpha -1)^2+\pi^2}\right) \, \eta(s) \, ds\\
&=&
\left[ \frac{\pi}{(s+t+\alpha -1)^2+\pi^2}\eta(s)\right] _0^{\infty}\\
& & - 
\int_0^{\infty}\frac{\pi}{(s+t+\alpha -1)^2+\pi^2}\, \eta'(s) \, ds\\
&=&
-\int_0^{\infty}\frac{\pi}{(s+t+\alpha -1)^2+\pi^2}\, \eta'(s) \, ds.
\end{eqnarray*}
If we replace $t$ by $-t$ in this formula we get
$$
Y_{-\pi}'(\alpha -t)=-\int_0^{\infty}\frac{\pi}{(s-t+\alpha -1)^2+\pi^2}\, \eta'(s) \,
ds,
$$
so that 
\begin{eqnarray*}
\lefteqn{Y_{-\pi}'(\alpha+t)-Y_{-\pi}'(\alpha-t)}\\
&=&
4\pi t\int_0^{\infty}\frac{s+\alpha -1}{((s+t+\alpha -1)^2+\pi^2)((s-t+\alpha -1)^2+\pi^2)}\, \eta'(s) \,
ds.
\end{eqnarray*}
Since $\eta'(s)>0$ we see from this expression that
$Y_{-\pi}'(\alpha+t)-Y_{-\pi}'(\alpha-t)$ has the same sign as $t$ for
$\alpha \geq 1$ and hence is
positive for $t>0$. \hfill $\square$

\section{Relation to work of Ramanujan}
\label{sec:ramanujan}
Defining $\theta(n)$ for natural
numbers $n$ by
\begin{equation}\label{eq:Ram}
  \frac{e^n}{2}=\sum_{k=0}^{n-1} \frac{n^k}{k!} +\theta(n)\frac{n^n}{n!},
\end{equation}
Ramanujan \cite{Ram} claimed that $\frac13<\theta(n)<\frac12$. This was later
proved independently by Szeg\H o \cite{Sz} and Watson
\cite{Wa}. Further details about Ramanujan's problem can be found in
\cite{B:C:K}. See also \cite{F:G:K:P}. Choi noticed
the relation between this problem and the median in the gamma distribution, see
\cite{choi}. By $n$ partial integrations we get the formula
$$
\frac{1}{n!}\int_0^n t^ne^{-t}\,dt=1-e^{-n}\sum_{k=0}^n\frac{n^k}{k!},
$$
hence by (\ref{eq:Ram})

$$
\frac{1}{n!}\int_0^n t^ne^{-t}\,dt=\frac{1}{2}-\frac{n^n}{n!}
(1-\theta(n))e^{-n},
$$
and finally
\begin{equation}\label{eq:Choi}
1-\theta(n)=\left(\frac{e}{n}\right)^n\int_n^{m(n+1)}t^ne^{-t}\,dt.
\end{equation}
Following Watson \cite{Wa} we can write
\begin{eqnarray*}
\theta(n)&=&1+\frac{e^n}{2n^n}n!-\frac{n!}{n^n}\sum_{k=0}^n \frac{n^k}{k!}\\
&=&1+\frac{e^n}{2n^n}\int_0^\infty t^ne^{-t}\,dt-\frac{1}{n^n}\int_0^\infty
(n+t)^ne^{-t}\,dt\\
&=&1+\frac{n}{2}\int_0^\infty t^ne^{n(1-t)}\,dt-\int_0^\infty (1+t/n)^n
e^{-t}\,dt.
\end{eqnarray*}

We now make the substitution $s=1+t/n$ in the last integral and split the
first integral in two by integrating over $(0,1)$ and $(1,\infty)$. This gives
\begin{equation}\label{eq:Wa}
  \theta(n)=1+\frac{n}{2}\left(\int_0^1 (te^{1-t})^n\,dt-
  \int_1^\infty (te^{1-t})^n\,dt\right).
\end{equation}
Watson notices that the right-hand side of this equation makes perfect
sence when $n$ is replaced by a positive real variable $x$ and the corresponding
function he denotes $y(x)$. We keep the notation $\theta(x)$ and make the
substitution $t=e^{-u}$ in the integrals to get
$$
\theta(x)=1+\frac{x}{2}e^x\left(\int_0^\infty e^{-x(u+e^{-u})-u}\,du 
-\int_{-\infty}^0 e^{-x(u+e^{-u})-u}\,du \right).
$$
Introducing the functions $u,v$ from Section \ref{sec:uniform} we get
$$
\theta(x)=1+\frac{x}{2}e^x\int_1^\infty e^{-xt}
\left(e^{-v(t)}v'(t)+e^{-u(t)}u'(t)\right)\,dt
$$
and using
$$
e^{-v(t)}v'(t)+e^{-u(t)}u'(t)=\xi(t)-2,
$$
we get after some calculation
\begin{equation}\label{eq:Wa2}
\theta(x)=\frac{x}{2}\int_0^\infty e^{-xt}\xi(t+1)\,dt=\frac13 +\frac12
\int_0^\infty e^{-xt}\xi'(t+1)\,dt.
\end{equation}
From (\ref{eq:partial_integration}) we get 
$$
\theta(x)=\sum_{k=0}^n\frac{\xi^{(k)}(1)}{2x^k} + o(x^{-n}).
$$
We further see that the function $\phi(t)$ from \cite[p.300]{Wa} is given
as
$$
\phi(t)=\frac{135}{8}\xi'(t+1).
$$
In \cite{Wa} Watson discusses the behaviour of $-t^{-1}\log \phi(t)$
for $t>0$. He reaches a conclusion using a combination of rigorous
analysis and numerical tabulation. The present analysis does not
contribute in making his conclusion rigorous.

\input maple_procedures.tex

\author{Department of Mathematics, University of Copenhagen,
  Universitetsparken 5, DK-2100, Copenhagen, Denmark.}\\
\noindent
\author{Email: berg@math.ku.dk\\Fax: +4535320704}
 \bigskip

\noindent
\author{Department of Natural Sciences, Royal Veterinary and
  Agricultural University,
  Thorvaldsensvej 40, DK-1871, Copenhagen, Denmark.}\\
\noindent
\author{Email: henrikp@dina.kvl.dk\\Fax: +4535282350}
\end{document}

%% file: maple_procedures.tex
\section{Appendix: Higher order expansions and Maple code}
In this section we give the derivatives of $\xi$ at 1 of order up to
10 and higher order asymptotic expansions of $\varphi(x)$ and $m(x)$. These results we have found using the ``Maple
9'' system. We have also included the Maple code, with a short
description.

\subsection*{Computation of the derivatives of $\xi$}
With the notation of Proposition \ref{prop:xi_at_1} we have
$\xi^{(k)}(1)=2a_{2(k+1)}(k+1)!$. Here the numbers $a_k$ are
defined by
$$
h^{-1}(w)=\sum_{k=1}^{\infty}a_{k}w^k,
$$
where $h^{-1}$ is the inverse to 
$$
h(z)=\sum_{k=1}^{\infty}b_{k}z^k,
$$
which satisfies $(h(z))^2=e^{-z}+z-1$ and $h'(0)>0$.

\bigskip

\noindent
The procedure \verb!h_polynomial(n)! returns the list of numbers
$[b_1, \ldots, b_n]$.
These numbers are computed by equating the
coefficients in the relation
$$
\left( \sum_{k=1}^{\infty}b_{k}z^k \right)^2 =
\sum_{k=2}^{\infty}\frac{(-1)^k}{k!}z^k.
$$
We get 
$$
\sum_{l=1}^{k-1}b_lb_{k-l}=\frac{(-1)^k}{k!}
$$
for $k\geq 2$, from which we determine the number $b_{k+1}$ from
$b_1,\ldots,b_k$. The numbers $b_1,\ldots,b_n$ are just the coefficients in the Taylor
expansion of $h(z)$ at $z=0$ up to order $n$. We could therefore also have used
the command \verb!taylor! but we prefer the more direct and efficient
method outlined in the procedure below.
\begin{verbatim}
h_polynomial:=proc(n)
        local b, j:
        b[1]:=1/sqrt(2):
        for j from 1 to (n-1) do
                b[j+1]:=((-1)^j/(j+2)!
                -sum(b[l+1]*b[j-l+1], l=1..(j-1)))/2/b[1]:
        end do: j:='j':
        return [seq(b[j], j=1..(n))];
end proc;
\end{verbatim}

\bigskip

\noindent
The procedure \verb!xi_derivatives(n)! computes the numbers
$[\xi(1),\ldots,\xi^{(n-1)}(1)]$ and it uses the procedure \verb!h_polynomial!. 
To find these numbers it is enough
to determine the numbers $a_1,\ldots,a_{2n}$ and this we do by
equating coefficients in the relation $h^{-1}(h(z))=z$. We put
$$
\sum_{l=k}^{\infty}b_{k,l}z^l= \left( \sum_{l=1}^{\infty}b_{l}z^l \right)^k. 
$$
Hence $b_{1,l}=b_l$ for $l\geq 1$ and by Cauchy multiplication,
$$
b_{k+1,l}=\sum_{m=k}^{l-1}b_{k,m}b_{l-m},
$$
for $l\geq k+1$.
Since $h^{-1}(h(z))=z$ we must furthermore have
$$
z = \sum_{k=1}^{\infty}a_k\left( \sum_{l=k}^{\infty}b_{k,l}z^l\right)
= \sum_{l=1}^{\infty}\left( \sum_{k=1}^la_kb_{k,l}\right) z^l. 
$$
Therefore $a_1b_{1,1}=1$ and 
$$
\sum_{k=1}^la_kb_{k,l}=0
$$
for $l\geq 2$. Hence $a_1=1/b_{1,1}=1/b_1$. 
The point is now that $a_l$ can be computed when we know
$a_1, \ldots ,a_{l-1}$ as well as $b_{1,l},\ldots ,b_{l,l}$. This is the main idea behind the procedure below.
\begin{verbatim}
xi_derivatives:=proc(N)
        local i,p,j,k,l,temp,a,h,b:
        h:=h_polynomial(2*N):
        for i from 1 to 2*N do
                for p from 1 to 2*N do
                        b[i,p]:=0:
                end do: p:='p':
        end do: i:='i':
        for j from 1 to 2*N do b[1,j]:=h[j]: end do: j:='j': 
        a[1]:=sqrt(2):
        for k from 1 to (2*N-1) do 
                for l from (k+1) to 2*N do
                        temp:=0:
                        for j from 1 to (l-1) do
                                temp:=temp+b[k,j]*h[l-j]:
                        end do: j:='j': 
                        b[k+1,l]:=temp:
                end do: l:='l':
                a[k+1]:=-sum(a[m]*b[m,k+1],m=1..k)/b[k+1,k+1]:
        end do: k:='k':
        return [seq(2*a[2*j]*(j!),j=1..N)];
end proc;
\end{verbatim}
We could also have used the command \verb!taylor! to compute the
coefficients in the expansion of $h^{-1}(h(z))$,
but it runs very slowly, so a more efficient program is needed.
We have found the derivatives of $\xi$ to be:
\begin{eqnarray*}
\xi(1) & = & \frac{2}{3}\\
\xi'(1) & = & \frac{2^3}{3^3\cdot 5}\\
\xi^{(2)}(1) & = & -\frac{2^4}{3^4\cdot 5\cdot 7}\\
\xi^{(3)}(1) & = & -\frac{2^5}{3^5\cdot 5\cdot 7}\\
\xi^{(4)}(1) & = & \frac{2^6\cdot 281}{3^8\cdot 5^2\cdot 7\cdot 11}\\
\xi^{(5)}(1) & = & \frac{2^7\cdot 23\cdot 227}{3^9\cdot 5^2\cdot
  7\cdot 11\cdot 13}\\
\xi^{(6)}(1) & = & -\frac{2^8\cdot 53\cdot 103}{3^{10}\cdot 5^2\cdot
  7\cdot 11\cdot 13}\\
\xi^{(7)}(1) & = & -\frac{2^9\cdot 373\cdot 4439\cdot 557}{3^{12}\cdot 5^4\cdot
  7^2\cdot 11\cdot 13\cdot 17}\\
\xi^{(8)}(1) & = & \frac{2^{10}\cdot 2650986803}{3^{13}\cdot 5^4\cdot
  7^2\cdot 11\cdot 13\cdot 17\cdot 19}\\
\xi^{(9)}(1) & = & \frac{2^{11}\cdot 6171801683}{3^{14}\cdot 5^4\cdot
  7^2\cdot 11\cdot 13\cdot 17\cdot 19}\\
\xi^{(10)}(1) & = & -\frac{2^{12}\cdot 1117\cdot 3835213201}{3^{16}\cdot 5^2\cdot
  7^2\cdot 11\cdot 13\cdot 17\cdot 19\cdot 23}.
\end{eqnarray*}

\subsection*{The asymptotic expansion of $\varphi$}
The asymptotic expansion of $\varphi(x)$ can be found from relation
(\ref{eq:phi_asymptotic}), that is
$$
\varphi(x)= 
\sum_{k=0}^{n}\frac{\xi^{(k)}(1)}{2x^{k+1}}- \sum_{k=1}^{n}\frac{(-1)^kx^k}{k!}\int_0^{\varphi(x)}(u+e^{-u}-1)^k\, du + o(x^{-(n+1)}).
$$
We use two procedures \verb!sum_of_int(n,PHI)! and
\verb!xi_sum(n)!. The procedure \verb!sum_of_int(n,PHI)! takes an expression of the form
$$
\verb!PHI! = \sum_{k=1}^{m}\frac{c_k}{x^{k}}
$$
and computes the sum of integrals
$$
\sum_{k=1}^{n}\frac{(-1)^kx^k}{k!}\int_0^{\verb!PHI!}(u+e^{-u}-1)^k\,
du 
$$
up to $o(x^{-(n+1)})$. This is done in the same way as in Lemma \ref{lemma:phi_asymptotic}. We put
$$
(e^{-u}+u-1)^k=\sum_{l=2k}^{\infty}t_{k,l}u^l,
$$
and observe that by Cauchy multiplication
$$
t_{k+1,m}=\sum_{l=2}^{m-2k}\frac{(-1)^l}{l!}t_{k,m-l},
$$
for $m\geq 2(k+1)$. Since $t_{1,m}=(-1)^m/m!$ we can compute the 
numbers $t_{k,\cdot}$ from the numbers $t_{k-1,\cdot}$ for $k\geq
2$. The numbers $t_{k,\cdot}$ are then used to approximate the
integrand $(u+e^{-u}-1)^k$ by its Taylor polynomial centered at $u=0$ of order $n+k$.
\begin{verbatim}
sum_of_int:=proc(n,PHI)
        local l,t,T,s,k,m,j,temp,temp1:
        for l from 2 to (n+1) do
                t[1,l]:=(-1)^l/l!:
        end do: l:='l':
        T:=sum(t[1,l]*PHI^(l+1)/(l+1), l=2..(n+1)):
        s:=-x*T:	
        for k from 2 to n do
                T:=0:
                for m from 2*k to (n+k) do
                        t[k,m]:=
                        sum((-1)^l/l!*t[k-1,m-l], l=2..(m-2*k+2)):
                        T:=T+t[k,m]*PHI^(m+1)/(m+1):
                end do: m:='m':
                s:=s+(-1)^k*x^k/k!*T:
        end do: k:='k':
        temp:=subs(x=1/x,expand(s)):
        temp1:=0:
        for j from 0 to (n+1) do
                temp1:=temp1+coeff(temp,x,j)*x^(-j):
        end do: j:='j':
        return temp1;	
end proc;
\end{verbatim}
We notice that the variable \verb!s! in the procedure
\verb!sum_of_int(n,PHI)! contains many terms of order $x^{-l}$ with
$l>k$. These terms are not needed and they are therefore
deleted. That happens in the \verb!for! loop containing \verb!temp1!.
(It improves program efficiency.)

\bigskip

\noindent
The second procedure \verb!xi_sum(n+1)! simply computes the asymptotic
expansion 
$$
\sum_{k=0}^{n}\frac{\xi^{(k)}(1)}{2x^{k+1}}
$$
by calling the procedure \verb!xi_derivatives(n+1)!.
\begin{verbatim}
xi_sum:=proc(n)
        local temp:
        temp:=xi_derivatives(n):
        return sum(temp[k]/2/x^k, k=1..n);
end proc;
\end{verbatim}

\bigskip

\noindent
To find the expansion of $\varphi(x)$ we combine these two procedures
in the procedure \verb!asymp_phi(n)!. This procedure computes the
asymptotic expansion of $\varphi(x)$ up to order $o(x^{-n-1})$.
\begin{verbatim}
asymp_phi:=proc(n)
        local phi,k:
        phi:=1/3/x:
        for k from 1 to n do
                phi:=xi_sum(k+1)-sum_of_int(k,phi):
        end do:
        return phi;
end proc;
\end{verbatim}
We have 
$$
\varphi(x)=\sum_{k=1}^{10}\frac{c_k}{x^{k}}+o(x^{-10}),
$$
where
\begin{eqnarray*}
c_1 & = & \frac{1}{3}\\
c_2 & = & \frac{29}{3^4\cdot 5 \cdot 2}\\
c_3 & = & -\frac{37}{3^6\cdot 5 \cdot 7}\\
c_4 & = & -\frac{3877}{3^9\cdot 5^2 \cdot 2^2}\\
c_5 & = & \frac{8957413}{3^{13}\cdot 5^3 \cdot 7\cdot 11}\\
c_6 & = & \frac{401\cdot 8842279}{2\cdot 3^{15}\cdot 5^2 \cdot
  7^2\cdot 11\cdot 13}\\
c_7 & = & -\frac{356146891\cdot 2039}{3^{18}\cdot 5^4 \cdot
  7^2\cdot 11\cdot 13}\\
c_8 & = & -\frac{216607304027\cdot 3077479}{2^3\cdot 3^{21}\cdot 5^6 \cdot
  7^3\cdot 11\cdot 13\cdot 17}\\
c_9 & = & \frac{31\cdot 743\cdot 4569027042343}{3^{23}\cdot 5^3 \cdot
  7^3\cdot 11\cdot 13\cdot 17\cdot 19}\\
c_{10} & = & \frac{71\cdot 282699240672481\cdot 1949 \cdot 5113}{2\cdot 3^{27}\cdot 5^7 \cdot
  7^3\cdot 11^2\cdot 13\cdot 17\cdot 19}.
\end{eqnarray*}

\subsection*{The asymptotic expansion of $m$}
Since $m(x)=xe^{-\varphi(x)}$ it is easy to find the asymptotic expansion of $m(x)$ when we have the
asymptotic expansion of $\varphi(x)$. It is done in the procedure
\verb!asymp_m(n)!, that gives the expansion of $m(x)$ up to order $x^{1-n}$.
\begin{verbatim}
asymp_m:=proc(n)
        local temp:
        temp:=subs(x=1/x, asymp_phi(n)):
        temp:=convert(taylor(exp(-temp), x=0, n+2), polynom):
        temp:=subs(x=1/x,temp):
        return expand(x*temp); 
end proc;
\end{verbatim}
We have 
$$
m(x)=x-\frac{1}{3}+\sum_{k=1}^9\frac{m_k}{x^{k}}+o(x^{-9}),
$$
where
\begin{eqnarray*}
m_1 & = & \frac{2^3}{3^4\cdot 5}\\
m_2 & = & \frac{2^3\cdot 23}{3^6\cdot 5 \cdot 7}\\
m_3 & = & \frac{2^3\cdot 281}{3^9\cdot 5^2 \cdot 7}\\
m_4 & = & -\frac{2^3\cdot 17\cdot 139753}{3^{13}\cdot 5^3 \cdot 7\cdot 11}\\
m_5 & = & -\frac{2^3\cdot 708494947}{3^{15}\cdot 5^3 \cdot 7^2\cdot
  11\cdot 13}\\
m_6 & = & \frac{2^3\cdot 140814348739}{3^{18}\cdot 5^4 \cdot 7^2\cdot
  11\cdot 13}\\
m_7 & = & \frac{2^3\cdot 7663181003289047}{3^{21}\cdot 5^6 \cdot 7^3\cdot
  11\cdot 13\cdot 17}\\
m_8 & = & -\frac{2^3\cdot 653\cdot 1359581\cdot 759929\cdot 3307}{3^{23}\cdot 5^6 \cdot 7^3\cdot
  11\cdot 13\cdot 17\cdot 19}\\
m_9 & = & -\frac{2^3\cdot 29\cdot 1376560394479059407}{3^{27}\cdot 5^7 \cdot 7^3\cdot
  11^2\cdot 17}.
\end{eqnarray*}